\input amstex
\documentstyle{amsppt}
\magnification=1200

\TagsOnRight
\NoRunningHeads
\NoBlackBoxes

\define\Y{\Bbb Y}
\define\R{\Bbb R}
\define\Z{\Bbb Z}
\define\C{\Bbb C}
\define\X{\frak X}
\define\al{\alpha}
\define\be{\beta}
\define\la{\lambda}

\def\Mz{M_{z,z'}}
\define\Mnz{M^{(n)}_{z,z'}}
\define\Mzxi{M_{z,z',\xi}}
\define\Pzxi{\Cal P_{z,z',\xi}}
\define\pitxi{\pi_{t,\xi}}
\define\Conf{\operatorname{Conf}}
\define\Res{\operatornamewithlimits{Res}}
\define\w{\tfrac{\xi}{\xi-1}}
\define\tht{\thetag}
\define\Ga{\Gamma}
\define\Rhat{\widehat R}
\define\Shat{\widehat S}

\define\M{\frak M}
\define\wt{\widetilde}

\topmatter
\title
Distributions on partitions, point processes, and the hypergeometric
kernel 
\endtitle
\author Alexei Borodin and Grigori Olshanski
\endauthor
\abstract
We study a 3--parametric family of stochastic point processes on the
one--di\-men\-sio\-nal lattice originated from a remarkable family of
representations of the infinite symmetric group. We prove that the
correlation functions of the processes are given by determinantal
formulas with a certain kernel. The kernel can be expressed through
the Gauss hypergeometric function; we call it the hypergeometric
kernel. 

In a scaling limit our processes approximate the processes describing
the decomposition of representations mentioned above into
irreducibles. As we showed before, see math.RT/9810015, the
correlation functions of these limit processes also have
determinantal form with so--called Whittaker kernel. We show that the
scaling limit of the hypergeometric kernel is the Whittaker kernel. 

The integral operator corresponding to the Whittaker kernel is an
integrable operator as defined by Its, Izergin, Korepin, and Slavnov.
We argue that the hypergeometric kernel can be considered as a kernel
defining a `discrete integrable operator'. 

We also show that the hypergeometric kernel degenerates for certain
values of parameters to the Christoffel--Darboux kernel for Meixner
orthogonal polynomials. This fact is parallel to the degeneration of
the Whittaker kernel to the Christoffel--Darboux kernel for Laguerre
polynomials. 
\endabstract

\toc
\widestnumber\head{\S7.}
\head \S0. Introduction \endhead
\head \S1. Distributions on partitions. The grand canonical ensemble
\endhead
\head \S2. Determinantal point processes \endhead
\head \S3. Calculation of the correlation functions. The
hypergeometric kernel \endhead 
\head \S4. Connection with Meixner polynomials \endhead
\head \S5. Scaling limit of the hypergeometric kernel: the Whittaker
kernel \endhead
\head \S6. Integrable operators \endhead
\head \S7. Appendix: some relations for the Gauss hypergeometric function
\endhead
\head{} References \endhead
\endtoc

\endtopmatter

\document
\head \S0. Introduction \endhead

Let $\Y_n$ be the set of Young diagrams with $n$ boxes and 
$\Y=\Y_0\sqcup\Y_1\sqcup\Y_2\sqcup\dots$ be the set of all Young
diagrams. In this paper we study a 
remarkable family of probability distributions on $\Y_n$,
$n=0,1,2,\dots$. 

The whole picture depends on 2 parameters $z$ and $z'$ which satisfy
certain conditions, see \S1. For each pair $(z,z')$ we have a
probability distribution on every $\Y_n$, $n=0,1,2,\dots$, we denote
it by $M_{z,z'}^{(n)}$. Its value on a Young diagram $\la\in \Y_n$
with Frobenius coordinates $(p_1,\dots,p_d\,|\,q_1,\dots,q_d)$ has
the form 
$$
\gathered
\Mnz(\la)=\frac{n!}{(zz')_n}\,(zz')^d\,\\ \times
\prod_{i=1}^d\frac{(z+1)_{p_i}(z'+1)_{p_i}(-z+1)_{q_i}(-z'+1)_{q_i}}
{p_i!p_i!q_i!q_i!}\,
{\det}^2\left[\frac1{p_i+q_j+1}\right]\,,
\endgathered\tag 0.1
$$
where $(a)_k$ stands for $a(a+1)\cdots(a+k-1)$. 

The distributions $\Mnz$ have a representation--theoretic meaning.
 Let $S_n$ be the symmetric group of degree $n$, $S(\infty)$ be
the union of the groups $S_n$, and for $\la\in\Y_n$, let $\chi^\la$
denote the irreducible character of $S_n$ corresponding to $\la$.
According to \cite{KOV}, there exists a central positive definite
function $\chi^{(z,z')}$ on $S(\infty)$ such that, for any $n$, its
restriction to $S_n$ is 
$$
\chi^{(z,z')}_n=\sum_{|\la|=n} \Mnz (\la)\frac{\chi^\la}{\chi^\la(e)}\,.
$$
Moreover, the unitary representation
$T^{(z,z')}$ corresponding to $\chi^{(z,z')}$ admits a
nice geometric description (at least for $z'=\bar z$), see \cite{KOV}.
This representation--theoretic aspect was the original motivation of
our study but in the present paper we do not discuss it (see
\cite{P.I}).  

Let us associate to a Young diagram
$\la=(p_1,\dots,p_d\,|\,q_1,\dots,q_d)$ a set of $2d$ points in
$\Z'=\Z+\frac 12$ as follows: 
$$
\la=(p_1,\dots,p_d\,|\,q_1,\dots,q_d)\,\mapsto\,
\{p_1+\tfrac12,\dots,p_d+\tfrac12,-q_1-\tfrac12,\dots,-q_d-\tfrac12\}. 
$$
Then every probability measure on $\Y_n$ provides a probability 
measure on the set of all point configurations in $\Z'$
with equal number of points in $\Z'_+=\Z'\cap \R_+$ and
$\Z'_-=\Z'\cap\R_-$ and such that the total sum of absolute values of
coordinates is equal to $n$. 

Next, having a distribution on each $\Y_n$, we can mix them using a
distribution on the set $\{0,1,2,\dots\}$ of indices $n$, then we get
a probability distribution on $\Y$. 

Thus, we get a probability measure on the set of all point
configurations in $\Z'$ with equal number of points in $\Z'_+$ and
$\Z'_-$. According to standard terminology, we can say that we
defined a point process on $\Z'$.  

Following a certain analogy with statistical physics, one can
call the resulting object of the mixing procedure the {\it grand
canonical ensemble}, see \cite{V}. 

For our special distributions (0.1) we choose the mixing distribution
to be the negative binomial distribution 
$$
\operatorname{Prob}\{n\}=(1-\xi)^t\,\frac{(t)_n}{n!}\, \xi^n, \quad
t=zz', \tag0.2 
$$
where $\xi\in (0,1)$ is an additional parameter. (This choice is
explained by willing to remove the factor $\frac{n!}{(t)_n}$ from the
RHS of (0.1).) We shall denote by 
$\Pzxi$ the point process on $\Z'$ thus obtained.

The main result of this paper is the explicit computation of the
correlation functions of $\Pzxi$. It turns out that they are given by
determinantal formulas 
$$
\rho_n(x_1,\dots,x_n)=\det[K(x_i,x_j)]_{i,j=1}^n
$$
with a certain kernel $K(x,y)$ on $\Z'$. This kernel can be expressed
through the Gauss hypergeometric function. We call it the {\it
hypergeometric kernel}.  

Due to the representation theoretic origin of our problem, the
distributions $\Mnz$ have a number of additional properties. In
particular, as $n\to\infty$, they converge to a probability measure
on a certain limit object called the `Thoma simplex', see \cite{KOO}
and \S5 below.  \footnote {This is a kind of dual object to the infinite
symmetric group, see \cite{T}, \cite{VK}, \cite{KV}. The limit
measure is, actually, a spectral measure for the decomposition of the
representation $T^{(z,z')}$ into irreducibles, see \cite{KOV}. }

In terms of point processes, this implies  that after an appropriate
scaling the point processes $\Pzxi$ will converge, as $\xi\nearrow1$,
 to a certain point process on $\Bbb R^*=\R\setminus\{0\}$
derived from  the limit measure on the Thoma simplex (we shall
give a rigorous proof of this result in our next paper).  

This limit process has been thoroughly studied in our previous papers
\cite{P.I} -- \cite{P.V}  \footnote{A survey of the results is given in
\cite{BO1}.}. Its correlation functions also have determinantal form 
with so--called Whittaker kernel, see \cite{P.IV}. In \S5 we show directly
that the scaling limit of the hypergeometric kernel is the Whittaker
kernel.  

The integral operator defined by the Whittaker kernel belongs to the
class of integrable operators as defined by Its, Izergin, Korepin,
Slavnov \cite{IIKS}. We show that the operator corresponding to
the hypergeometric kernel can be considered as an example of a
`discrete integrable operator'.

A.~Okounkov pointed out that important information can be 
obtained from consideration of 
another degeneration of the point process introduced above. Assume that
$z,z'\to\infty$ and $\xi=\frac {\eta}{zz'}=\frac\eta t\to 0$ where
$\eta>0$ is fixed.
Then the mixing distribution (0.2) tends to the Poisson distribution with
parameter $\eta$:
$$
\operatorname{Prob} \{n\}\to e^{-\eta}\,\frac {\eta^n}{n!},
$$                                  
while $M_{zz'}^{(n)}$ tends to the {\it Plancherel distribution} on $\Bbb
Y_n$:
$$
M_{zz'}^{(n)}(\lambda)\to
M_{\infty}^{(n)}(\lambda)=\frac{\dim^2\lambda}{n!}\,, \tag0.3
$$
where 
$\dim\lambda=\chi^\lambda(e)$ is the dimension of the irreducible
representation of the symmetric group $S_n$ corresponding to
$\lambda$. Thus, we get an explicit formula for the
correlation functions of the process governed by the  poissonized
Plancherel distributions. \footnote{The limit relation \thetag{0.3}
was known since the invention of the distributions $M_{zz'}^{(n)}$,
see \cite{KOV}, but up to now it was not used.}

This formula allows to prove certain important facts
about Plancherel distributions, see \cite{BOO}.
In particular, we were able to prove the conjecture by Baik, Deift, and
Johansson \cite{BDJ1, BDJ2} that the asymptotic behavior of 
$\lambda_1,\lambda_2,\dots$ with respect to the Plancherel distribution is
governed by the Airy kernel \cite{TW} and, therefore, coincides with that of
the largest eigenvalues of a matrix from the Gaussian Unitary Ensemble.
(Another approach to this conjecture can be found in \cite{O}.)

The paper is organized as follows.
In \S1 we introduce our main object of interest -- the point process
$\Pzxi$. In \S2 we recall some generalities on determinantal point
processes. The computation of the correlation functions of $\Pzxi$
and the formulas for the hypergeometric kernel can be found in \S3.
In \S4 we show that if one of the parameters $z$, $z'$ is an integer,
the hypergeometric kernel degenerates to the Christoffel--Darboux
kernel for the Meixner orthogonal polynomials on $\Z_+$.
\footnote{After the present paper was completed we learned 
that the `Meixner kernel' has also arisen in the recent work
\cite{J}.} In \S5 we discuss the scaling limit of the hypergeometric
kernel. \S6 explains the connection with integrable operators. \S7 is an
appendix, we give there proofs of certain identities involving the Gauss
hypergeometric functions which are used in \S3. 

\subhead Acknowledgements \endsubhead We are grateful to P.~Deift and
N.~A.~Slavnov for consultations concerning the subject of \S6. We
also thank C.~A.~Tracy for telling us about Johansson's talk at MSRI
and sending us a copy of the transparencies, and K.~Johansson for
further information about his work \cite{J}.

\head \S1. Distributions on partitions. The grand canonical ensemble
\endhead 

For $n=1,2,\dots$, let $\Y_n$ denote the set of partitions of $n$,
which will be identified with Young diagrams with $n$ boxes. We agree
that $\Y_0$ consists of a single element --- the zero partition or
the empty diagram $\varnothing$. 

Given $\la\in\Y_n$, we write $|\la|=n$ and denote by $d=d(\la)$ the
number of diagonal boxes in $\la$. We shall use the Frobenius
notation \cite{Ma} 
$$
\la=(p_1,\dots,p_d\,|\,q_1,\dots,q_d).
$$
Here $p_i=\la_i-i$ is the number of boxes in the $i$th row of $\la$ on the
right of the $i$th diagonal box; likewise, $q_i=\la'_i-i$ is the number of 
boxes in the $i$th column of $\la$ below the $i$th diagonal box ($\la'$
stands for the transposed diagram). 


Note that 
$$
p_1>\dots>p_d\ge0, \qquad
q_1>\dots>q_d\ge0, \qquad
\sum_{i=1}^d(p_i+q_i+1)=|\la|.
$$
The numbers $p_i$, $q_i$ are called the {\it Frobenius coordinates\/}
of the diagram $\la$. 

Throughout the paper we fix two complex parameters $z,z'$ such that
the numbers $(z)_k(z')_k$ and $(-z)_k(-z')_k$ are real and strictly
positive for any $k=1,2,\dots$. Here and below 
$$
(a)_k=a(a+1)\dots(a+k-1), \qquad (a)_0=1, 
$$ 
denotes the Pochhammer symbol. 

The above assumption on $z,z'$ means that one of the
following two conditions holds:

$\bullet$ either $z'=\bar z$ and $z\in\C\setminus\Z$

$\bullet$ or $z,z'\in\R$ and there exists $m\in\Z$ such that
$m<z,z'<m+1$. 

We set 
$$
t=zz'
$$
and note that $t>0$.

For a Young diagram $\la$ let $\dim\la$ denote the number of the
standard Young tableaux of shape $\la$. Equivalently, $\dim\la$ is
the dimension of the irreducible representation (of the symmetric
group of degree $|\la|$) corresponding to $\la$,  see
\cite{Ma}. In the Frobenius notation, 
$$
\frac{\dim\la}{|\la|!}=
\frac{\prod\limits_{1\le i<j\le d}(p_i-p_j)(q_i-q_j)}
{\prod\limits_{1\le i,j\le n}(p_i+q_j+1)}
=\det\left[\frac1{p_i+q_j+1}\right]_{1\le i,j\le d}\,,
$$
see, e.g., \cite{P.I, Prop. 2.6}. 

We introduce a function on the Young diagrams depending on the
parameters $z,z'$:
$$
\gather
\Mz(\la)=\frac{t^d}{(t)_n}\,
\prod_{i=1}^d\frac{(z+1)_{p_i}(z'+1)_{p_i}(-z+1)_{q_i}(-z'+1)_{q_i}}
{p_i!p_i!q_i!q_i!}\,
\frac{\dim^2\la}{|\la|!}\\
=\frac{|\la|!}{(t)_n}\,t^d\,
\prod_{i=1}^d\frac{(z+1)_{p_i}(z'+1)_{p_i}(-z+1)_{q_i}(-z'+1)_{q_i}}
{p_i!p_i!q_i!q_i!}\,
{\det}^2\left[\frac1{p_i+q_j+1}\right]\,.
\endgather
$$

We agree that $\Mz(\varnothing)=1$. Thanks to our assumption on the
parameters, $\Mz(\la)>0$ for all $\la$. 

\proclaim{Proposition 1.1} For any $n$,
$$
\sum_{\la\in\Y_n}\Mz(\la)=1,
$$
so that the restriction of $\Mz$ to $\Y_n$ is a probability
distribution on $\Y_n$.
\endproclaim

We shall denote this distribution by $\Mnz$.

\demo{Comments} This result is the starting point of our
investigations. About its origin and representation--theoretic
significance, see \cite{KOV}. Several direct proofs of the
proposition are known. E.g., a simple proof is given in \cite{Part I,
\S7}. About generalizations, see \cite{K}, \cite{BO2}. Note that 
$$
\lim_{|z|,|z'|\to\infty}\Mz(\la)=\frac{\dim^2\la}{|\la|!}\,,
$$
so that the limit form of the identity of the proposition is
$$
\sum_{\la\in\Y_n}\frac{\dim^2\la}{|\la|!}=1,
$$
which is well known.
\enddemo

Let $\Y=\Y_0\sqcup\Y_1\sqcup\dots$ denote the set of all Young
diagrams. Consider the negative binomial distribution on the
nonnegative integers, which depends on $t$ and the additional parameter
$\xi$, $0<\xi<1$:
$$
\pitxi(n)=(1-\xi)^t\,\frac{(t)_n}{n!}\,\xi^n\,,
\qquad n=0,1,\dots\,.
$$
For $\la\in\Y$ we set
$$
\Mzxi(\la)=\Mz(\la)\,\pitxi(|\la|).
$$
By the construction, $\Mzxi(\cdot)$ is a probability distribution on
$\Y$, which can be viewed as a mixing of the finite distributions
$\Mnz$. {}From the formulas for $\Mz$ and $\pitxi$ we get an explicit
expression for $\Mzxi$:
$$
\multline
\Mzxi(\la)=(1-\xi)^t\,
{\xi\,}^{\sum\limits_{i=1}^d(p_i+q_i+1)}\,
t^d\,\\
\times\prod_{i=1}^d\frac{(z+1)_{p_i}(z'+1)_{p_i}(-z+1)_{q_i}(-z'+1)_{q_i}}
{p_i!p_i!q_i!q_i!}\,
{\det}^2\left[\frac1{p_i+q_j+1}\right]\,.
\endmultline \tag1.1
$$

Following a certain analogy with models of statistical physics (cf.
\cite{V}) one may call $(\Y,\Mzxi)$ the {\it grand canonical
ensemble.}

Let $\Z'$ denote the set of half--integers,
$$
\Z'=\Z+\tfrac12=\{\dots,-\tfrac32,-\tfrac12,\tfrac12,\tfrac32,\dots\},
$$
and let $\Z'_+$ and $\Z'_-$ be the subsets of positive and negative
half--integers, respectively.
It will be sometimes convenient to identify both
$\Z'_+$ and $\Z'_-$ with $\Z_+=\{0,1,2,\dots\}$ by making use of the
correspondence $\pm(k+\tfrac12)\leftrightarrow k$, where $k\in\Z_+$.

Denote by $\Conf(\Z')$ the space of all finite subsets of $\Z'$ which
will be called {\it configurations}. We define an embedding
$\la\mapsto X$ of the set $\Y$ of Young 
diagrams into the set $\Conf(\Z')$ of configurations in $\Z'$ as follows: 
$$
\la=(p_1,\dots,p_d\,|\,q_1,\dots,q_d)\,\mapsto\,
X=\{p_1+\tfrac12,\dots,p_d+\tfrac12,-q_1-\tfrac12,\dots,-q_d-\tfrac12\}. 
$$
Under the identification $\Z'\simeq \Z_+\sqcup \Z_+$, the map
$\la\mapsto X$ is simply associating to $\la$ the collection of its
Frobenius coordinates. The image of the map consists exactly of the
configurations $X$ with the property $|X\cap\Z'_+ |=|X\cap \Z'_-|$.

 Under the embedding $\la\mapsto X$ the probability measure $\Mzxi$ on
$\Y$ turns into a probability
measure on the configurations in $\Z'$. Following the conventional 
terminology, see \cite{DVJ}, we get a point process on $\Z'$; let us
denote it as $\Pzxi$. 

Our primary goal is to compute the correlation functions of this
point process. 

\head \S2. Determinantal point processes \endhead

Let $\X$ be a countable set. Its finite subsets will be called {\it
configurations\/} and denoted by the letters $X,Y$. The space of all
configurations is denoted as $\Conf(\X)$; this is a discrete space.
In this section,  by a {\it point process\/} on $\X$ we mean a map
from a probability space to $\Conf(\X)$. \footnote{Actually, such a 
definition is rather restricted but it suffices for our purpose. For
general axiomatics of point processes, see \cite{DVJ}, \cite{L1},
\cite{L2}.} Let $\Cal P$ be a point process on $\X$. It induces a
probability distribution on $\Conf(\X)$, which is a nonnegative function
$\pi(X)$ such that $\sum_X\pi(X)=1$. One may simply identify $\Cal P$ and
$\pi$: then the underlying probability space is $\Conf(\X)$ itself.

We introduce a related function
on $\Conf(\X)$ as follows: 
$$
\rho(X)=\sum_{Y\supseteq X}\pi(Y).
$$
That is, $\rho(X)$ is the probability that the random configuration
contains $X$. 

Consider the Hilbert space $\ell^2(\X)$. An operator $L$ in
$\ell^2(\X)$ will be viewed as an infinite matrix $L(x,y)$ whose rows
and columns are indexed by points of $\X$. By $L_X$ we denote the
finite matrix of format $X\times X$ which is obtained from $L$ by
letting $x,y$ range over $X$. 

Assume $L$ is a trace class operator in $\X$ such that all the
principal minors $\det L_X$ are real nonnegative numbers. We agree
that $\det L_\varnothing=1$. We have
$$
\det(1+L)=\sum_n\operatorname{tr}(\wedge^n L)=
\sum_X\det L_X
$$
(here $\wedge^n L$ stands for the $n$th exterior power of $L$ acting
in the $n$th exterior power of the Hilbert space $\ell^2(\X)$). 
Thus, we can define a point process by
$$
\pi(X)=\frac{\det L_X}{\det(1+L)}\,.
$$
Let us call it the {\it determinantal point process\/} determined by
the operator $L$. \footnote{This term is not a conventional one. Such
processes, not necessarily on discrete spaces, arise in different
topics, in particular, in connection with random matrices. 
In \cite{DVJ}, they are called `fermion processes' but in
the the random matrix literature no special term is adopted. }

\proclaim{Proposition 2.1} Let $L$ satisfy the above
assumption, $\pi(\cdot)$ be the corresponding point process, and
$\rho(\cdot)$ be the associated function as defined above. Set
$K=L(1+L)^{-1}$. Then 
$$
\rho(X)=\det K_X\,.
$$
\endproclaim

\demo{Proof} We shall reproduce the argument indicated in \cite{DVJ}.
 Let $f(x)$ be a function on $\X$ such that 
$f_0(x)=f(x)-1$ is finitely supported. For any point process
$\pi(\cdot)$, 
$$
\sum_Y\pi(Y)\prod_{y\in Y}f(y)=
\sum_Y\pi(Y)\prod_{y\in Y}(1+f_0(y))=
\sum_X\rho(X)\prod_{x\in X}f_0(x).
$$
When the process is defined by an operator $L$ then, identifying $f$
with the diagonal matrix with diagonal entries $f(x)$, we get
$$
\gather
\sum_Y\pi(Y)\prod_{y\in Y}f(y)
=\left(\sum_Y\det L_Y \prod_{y\in Y}f(y)\right)
{\det}^{-1}(1+L)\\
=\det(1+fL){\det}^{-1}(1+L)=
\det\left((1+fL)(1+L)^{-1}\right)\\
=\det\left((1+L+f_0L)(1+L)^{-1}\right)=
\det(1+f_0K)=
\sum_X\det K_X\prod_{x\in X}f_0(x)\,.
\endgather
$$
Thus,
$$
\sum_X\rho(X)\prod_{x\in X}f_0(x)=\sum_X\det K_X\prod_{x\in X}f_0(x)
$$
for any finitely supported function $f_0$, which implies
$\rho(X)=\det K_X$. \qed
\enddemo

Let $\rho_n$ be the restriction of $\rho(\cdot)$ to the $n$-point
configurations. One can view $\rho_n$ as a symmetric function in $n$
variables,
$$
\rho_n(x_1,\dots,x_n)=\rho(\{x_1,\dots,x_n\}),
\quad\text{$x_1,\dots,x_n$ pairwise distinct.}
$$
In this notation, the result of Proposition 2.1 reads as follows
$$
\rho_n(x_1,\dots,x_n)=\det[K(x_i,x_j)]_{1\le i,j\le n}\,.
$$
We call $\rho_n$ the {\it $n$-point correlation function.}

{}From now on we assume that $\X=\X^+\sqcup\X^-$ (disjoint union of two
countable sets) and we write
$\ell^2(\X)=\ell^2(\X^+)\oplus\ell^2(\X^-)$. According to this
decomposition we write operators in $\ell^2(\X)$ as $2\times2$
operator matrices,
$$
L=\bmatrix L_{++} & L_{+-}\\ L_{-+} & L_{--}
\endbmatrix\,, \qquad
K=\bmatrix K_{++} & K_{+-}\\ K_{-+} & K_{--}
\endbmatrix\,.
$$

Given a configuration $X$, we set $X^\pm=X\cap\X^\pm$. We shall deal
with operators $L$ such that $L_{++}=0$, $L_{--}=0$. Then, as is
easily seen, $\pi(X)=0$ unless $|X^+|=|X^-|$. 

\proclaim{Proposition 2.2} The transforms $L\mapsto K=L(1+L)^{-1}$
and $K\mapsto L=K(1-K)^{-1}$ 
define a bijective correspondence between

{\rm(i)} the operators $L$ of the form
$$
L=\bmatrix 0 & A\\ -B & 0\endbmatrix\,, \tag2.1
$$
where the matrix $1+AB$ is invertible (equivalently, $1+BA$ is
invertible) 

\noindent and

{\rm(ii)} the operators $K$ of the form
$$
K=\bmatrix CD & C\\ DCD-D & DC\endbmatrix \,, \tag2.2
$$
where $1-CD$ is invertible (equivalently, $1-DC$ is invertible).

In terms of the blocks, this correspondence takes the form
$$
\gather
C=(1+AB)^{-1}A=A(1+BA)^{-1}, \quad D=B,\tag2.3\\
A=C(1-DC)^{-1}=(1-CD)^{-1}C , \quad B=D.\tag2.4
\endgather
$$
In particular, 
$$
1-CD=(1+AB)^{-1}, \quad 1-DC=(1+BA)^{-1}.\tag2.5
$$
\endproclaim

\demo{Proof} The proof is straightforward, see \cite{P.V, Prop.
2.2}.  \qed 
\enddemo

\proclaim{Proposition 2.3} Let $K$ be a $J$--Hermitian
\footnote{I.e., Hermitian with respect to the indefinite inner
product determined by the matrix 
$J=\bmatrix 1 & 0\\ 0 & -1\endbmatrix$. }
kernel of the
form \tht{2.2}. Then 
$$
L=\bmatrix 0 & D^* \\ -D & 0 \endbmatrix\,.
$$
\endproclaim

\demo{Proof} By Proposition 2.2, $L$ is given by the formula
\tht{2.1} with $B=D$. Since $K$ is $J$-Hermitian, $L$ is
$J$-Hermitian, too. This implies $A=B^*=D^*$. \qed
\enddemo

\head \S3. Calculation of the correlation functions. The
hypergeometric kernel \endhead 

In this section we shall apply the
formalism of \S2 to the point processes $\Pzxi$ introduced at the end of \S1.

Naturally, we specify $\X=\Z'$ and $\X^\pm=\Z'_\pm$. 

Let us introduce 2 meromoprhic functions in $u$ depending on $z$,
$z'$, $\xi$ as parameters 
$$
\gathered
\psi_\pm(u)=t^{1/2}\,\xi^{u+1/2}\,(1-\xi)^{\pm(z+z')}\,
\frac{\Ga(u+1\pm z)\Ga(u+1\pm z')}
{\Ga(1\pm z)\Ga(1\pm z')\Ga(u+1)\Ga(u+1)}   \\
=t^{-1/2}\,\xi^{u+1/2}\,(1-\xi)^{\pm(z+z')}\,
\frac{\Ga(u+1\pm z)\Ga(u+1\pm z')}
{\Ga(\pm z)\Ga(\pm z')\Ga(u+1)\Ga(u+1)}.
\endgathered
\tag3.1 
$$

An important fact is that the functions $\psi_\pm(u)$ have
exponential decay as $u$ tends to $+\infty$ along the real axis.
(Indeed, since $\xi\in(0,1)$, the factor $\xi^u$ has exponential 
decay; as for the remaining expression in \tht{3.1}, it behaves as a
constant times $u^{\pm(z+z')}$, so that it has at most polynomial
growth.)

We shall consider two diagonal matrices of format $\Z_+\times\Z_+$,
depending on the parameters $z,z',\xi$ and denoted as
$\Psi_\pm$, and a third matrix of the same format denoted as $W$:
$$
\Psi_\pm=\operatorname{diag}
\left\{\frac{t^{1/2}\xi^{k+1/2}(1-\xi)^{\pm(z+z')}
(1\pm z)_k(1\pm z')_k}{k!\, k!}\right\}\,, \quad
W(k,l)=\frac1{k+l+1}\,, \tag3.2
$$
where $k,l$ range over $\Z_+$. 
Note that the $k$th diagonal entry of $\Psi_\pm$ equals
$\psi_\pm(k)$. As the diagonal entries of $\Psi_\pm$ are real
and positive, we may introduce the diagonal matrices $\Psi^{1/2}_\pm$
which are real and positive, too. Note also that $W$ is real and
symmetric.

\proclaim{Proposition 3.1} The point process $\Pzxi$ is
a determinantal process in the sense of {\rm\S2}, and the
corresponding operator $L$ is as follows: 
$$
L=\bmatrix 0 & \Psi^{1/2}_+ W \Psi^{1/2}_- \\
-\Psi^{1/2}_- W \Psi^{1/2}_+ & 0 \endbmatrix \,. \tag3.3
$$
\endproclaim

Note that $L$ is real and $J$--symmetric. 

\demo{Proof} Let $\la$ be a Young diagram and $X$ be the
corresponding configuration. We must prove that
$$
\Mzxi(\la)=\frac{\det L_X}{\det(1+L)}\,. \tag3.4
$$
Since $L$ has the form $\bmatrix 0 & A \\ -A' & 0 \endbmatrix\,$,
where prime means transposition, we have, taking account of the exact
expression for the matrix $A$ (see \tht{3.3}):
$$
\det L_X={\det}^2 [A(p_i,q_j)]_{1\le i,j\le d}
=\prod_{i=1}^d\psi_+(p_i)\psi_-(q_i)\cdot
{\det}^2\left[\frac1{p_i+q_j+1}\right]\,. 
$$
In the latter product, the factor $(1-\xi)^{z+z'}$ coming from
$\psi_+$ cancels with the factor $(1-\xi)^{-(z+z')}$, and we get
$$
\multline
\det L_X=\det(1+L)^{-1}\,
{\xi\,}^{\sum\limits_{i=1}^d(p_i+q_i+1)}\,
t^d\,\\
\times\prod_{i=1}^d\frac{(z+1)_{p_i}(z'+1)_{p_i}(-z+1)_{q_i}(-z'+1)_{q_i}}
{p_i!p_i!q_i!q_i!}\,
{\det}^2\left[\frac1{p_i+q_j+1}\right]\,.
\endmultline
$$
Comparing this with the expression \tht{1.1} for $\Mzxi(\la)$, we see
that they coincide up to a constant factor which does not depend on $\la$.
Since the both expressions define probability distributions, we
conclude that they are identical.  \qed
\enddemo

\example{Remark 3.2} As a by--product of the proof we get the
following result: 
$$
\det(1+L)=\det(1+\Psi^{1/2}_+W\Psi_-W\Psi^{1/2}_+)
=(1-\xi)^{-t}\,.
$$
\endexample

By Propositions 3.1 and 2.1, the correlation functions of the point
process $\Pzxi$ are given by the determinantal formula involving
the operator $K=L(1+L)^{-1}$. Theorem 3.3 below provides an explicit
expression for this operator. 

Introduce the functions 
$$
R_\pm(u)=\psi_\pm(u)\, F(\mp z,\mp z';u+1;\w), 
\tag 3.5
$$
$$
S_\pm(u)=\frac{t^{1/2}\xi^{1/2}\,\psi_\pm(u)}{1-\xi}\,
\frac{F(1\mp z,1\mp z';u+2;\w)}{u+1},
\tag 3.6
$$
$$
P_\pm(u)=(\psi_\pm(u))^{-1/2}R_\pm(u),\quad
Q_\pm(u)=(\psi_\pm(u))^{-1/2}S_\pm(u). 
\tag 3.7
$$
They are all well--defined for $u\ge0$, because $\psi_\pm(u)$ is strictly
positive for $u\ge0$. In particular, they are well--defined at the
points $u=k\in\Z_+$. 

Note that the hypergeometric function that enters \tht{3.5} or
\tht{3.6} remains bounded as $u\to+\infty$. Hence, the exponential
decay of $\psi_\pm(u)$ implies the exponential decay of $R_\pm(u)$
and $S_\pm(u)$ as $u$ tends to $+\infty$ along the real axis.  

\proclaim{Theorem 3.3} Let 
$K=\bmatrix K_{++} & K_{+-} \\ K_{-+} & K_{--} \endbmatrix$ 
be the operator in $\ell^2(\Z_+)\oplus\ell^2(\Z_+)$ with the blocks
$$
K_{++}(k,l)=\frac{P_+(k)Q_+(l)-Q_+(k)P_+(l)}{k-l}\,, $$$$
K_{+-}(k,l)=\frac{P_+(k)P_-(l)+Q_+(k)Q_-(l)}{k+l+1}\,, $$$$
K_{-+}(k,l)=-\,\frac{P_-(k)P_+(l)+Q_-(k)Q_+(l)}{k+l+1}\,,$$$$
K_{--}(k,l)=\frac{P_-(k)Q_-(l)-Q_-(k)P_-(l)}{k-l}\,.
$$
Here the functions $P_\pm(u)$ and $Q_\pm(u)$ are defined in \tht{3.7},
and the expressions $K_\pm(k,l)\bigm|_{k=l}$ are understood according
to the L'Hospital rule.

Then $K=L(1+L)^{-1}$, where $L$ is defined in \tht{3.3}.
\endproclaim

We shall call the kernel $K$ defined above the {\it hypergeometric kernel}. 

Note that $K$ is $J$--symmetric (because of the minus sign in
the expression for $K_{-+}$). 

\demo{Proof} We shall prove that $K$ has the form \tht{2.2} with
$$
C=K_{+-}\,,  \qquad D=-L_{-+}=\Psi^{1/2}_-W\Psi^{1/2}_+ \,.  \tag 3.8
$$
I.e., 
$$
K_{++}=CD, \qquad K_{--}=DC, \qquad 
K_{-+}=DCD-D.
\tag 3.9
$$
As $K$ is $J$--symmetric, the desired result will follow from
Proposition 2.3. \footnote{In the same way, one could verify directly
that the operators $C,D$ obey the relations \tht{2.3} which means that
$K$ coincides with $L(1+L)^{-1}$. However, reference to Proposition
2.3 makes this verification redundant.}

It is convenient to slightly rewrite the desired relations \tht{3.9} in
order to avoid square roots. To do this, we set
$$
N=\Psi^{1/2}_+C\Psi^{1/2}_-=\Psi_+W\Psi_-\,.
$$
I.e., 
$$
N(k,l)=\frac{R_+(k)R_-(l)+S_+(k)S_-(l)}{k+l+1}\,.
$$
By virtue of the connection between $P_\pm, Q_\pm$ and $R_\pm,
S_\pm$, the relations \tht{3.9} are equivalent to the following ones:
$$
\gather
(NW)(k,l)=\frac1{\psi_+(l)}\,
\frac{R_+(k)S_+(l)-S_+(k)R_+(l)}{k-l}\,,   \tag3.10 \\
(WN)(k,l)=\frac1{\psi_-(k)}\,
\frac{R_-(k)S_-(l)-S_-(k)R_-(l)}{k-l}\,,   \tag3.11 \\
(WNW-W)(k,l)=\frac1{\psi_-(k)\psi_+(l)}\,
\frac{R_-(k)R_+(l)-S_-(k)S_+(l)}{k+l+1}\,.   \tag3.12 \\
\endgather
$$
We note once more that, by agreement, the indeterminacy arising in
\tht{3.10} and 
\tht{3.11} for $k=l$ is removed by making use of the L'Hospital rule.
 
To prove the relations above we shall need certain identities
involving the hypergeometric function.  

\proclaim{Lemma 3.4}
Set 
$$
\Rhat_\pm(u)=\sum_{k=0}^\infty \frac{R_\pm(k)}{u+k+1}\,, \qquad
\Shat_\pm(u)=\sum_{k=0}^\infty \frac{S_\pm(k)}{u+k+1}\,.  
\tag 3.13
$$
Then the series absolutely converge for $u\ne -1,-2,\dots$ and the
following relations hold
$$
\Rhat_\pm(u)=\psi_\mp(u)^{-1}\,S_\mp(u), \qquad
\Shat_\pm(u)=1-\psi_\mp(u)^{-1}\,R_\mp(u).   \tag3.14
$$
\endproclaim
\proclaim{Lemma 3.5}
$$
R_+(u)R_-(-u-1)+S_+(u)S_-(-u-1)= \psi_+(u)\psi_-(-u-1).\tag3.15
$$
\endproclaim
The proofs of these two lemmas can be found in the Appendix (\S7). 

Let us now  check \tht{3.10}. By the definition of $N$ and $W$,
$$
(NW)(k,l)=\sum_{j=0}^\infty 
\frac{R_+(k)R_-(j)+S_+(k)S_-(j)}{(k+j+1)(j+l+1)} \,.  
\tag 3.16
$$
Assume first $k\ne l$. We write
$$
\frac1{(k+j+1)(j+l+1)}=
\frac1{k-l}\left(\frac1{l+j+1}-\frac1{k+j+1}\right)\,,   
\tag 3.17
$$ 
and plug this into \tht{3.16}. Then we get
$$
\gather
(NW)(k,l)=\frac{R_+(k)}{k-l}\,\sum_{j=0}^\infty
\frac{R_-(j)}{l+j+1}\, + \,
\frac{S_+(k)}{k-l}\,\sum_{j=0}^\infty
\frac{S_-(j)}{l+j+1}\\
-\frac{R_+(k)}{k-l}\,\sum_{j=0}^\infty
\frac{R_-(j)}{k+j+1}\,-\,
\frac{S_+(k)}{k-l}\,\sum_{j=0}^\infty
\frac{S_-(j)}{k+j+1}\,.
\endgather
$$
By \tht{3.13}, this can be written as  
$$
(NW)(k,l)=
\frac{R_+(k)\Rhat_-(l)+S_+(k)\Shat_-(l)}{k-l}\,-\,
\frac{R_+(k)\Rhat_-(k)+S_+(k)\Shat_-(k)}{k-l}\,.
$$
Applying \tht{3.14}, we get
$$
\gather
R_+(k)\Rhat_-(l)+S_+(k)\Shat_-(l)
=\frac{R_+(k)S_+(l)-S_+(k)R_+(l)}{\psi_+(l)}\,,\\
R_+(k)\Rhat_-(k)+S_+(k)\Shat_-(k)
=\frac{R_+(k)S_+(k)-S_+(k)R_+(k)}{\psi_+(k)}=0.
\endgather
$$
Thus, we have checked \tht{3.10} for $k\ne l$. 

To extend the argument above to the case $k=l$, we replace \tht{3.17} by 
a slightly more complicated expression that makes sense for any
$k,l\in\Z_+$: 
$$
\frac1{(k+j+1)(j+l+1)}=\lim_{u\to l}
\left\{\frac1{k-u}\left(\frac1{u+j+1}-\frac1{k+j+1}\right)\right\}\,,
\tag 3.18
$$
where $u$ is assumed to be nonintegral. 
Since the functions $R_\pm(u)$, $S_\pm(u)$ have exponential decay as
$u\to+\infty$ (see the paragraph before Theorem 3.3), we may
interchange summation and the limit 
transition. Then we can repeat all the transformations. At the very
end we must pass to the limit as $u\to l$, which means that we follow
the L'Hospital rule. This concludes the proof of \tht{3.16}. 

The proof of \tht{3.11} is quite similar, and we proceed to the proof
of \tht{3.12}.  

By virtue of the expression \tht{3.10} for $NW$ and our agreement about the
L'Hospi\-tal rule, we get
$$
(WNW)(k,l)=\lim_{u\to l}\left\{\sum_{j=0}^\infty
\frac{R_+(j)S_+(u)-S_+(j)R_+(u)}{\psi_+(u)(k+j+1)(j-u)}\right\}\,. 
\tag 3.19
$$
Using the transformation
$$
\frac1{(k+j+1)(j-u)}=-\,\frac1{k+u+1}\,
\left(\frac1{k+j+1}-\,\frac1{(-u-1)+j+1}\right)
$$
we rewrite the above sum as follows:
$$
\multline
\sum_{j=0}^\infty
\frac{R_+(j)S_+(u)-S_+(j)R_+(u)}{\psi_+(u)(k+j+1)(j-u)}\\
=-\, \frac{\Rhat_+(k)S_+(u)-\Shat_+(k)R_+(u)}{k+u+1}
+\, \frac{\Rhat_+(-u-1)S_+(u)-\Shat_+(-u-1)R_+(u)}{k+u+1}\,.
\endmultline
$$
Next, applying \tht{3.14}, we transform this to 
$$
\multline
-\,\frac{S_-(k)S_+(u)+R_-(k)R_+(u)}{\psi_-(k)(k+u+1)}\,-\,
\frac{R_+(u)}{k+u+1}\\
+\, \frac{S_-(-u-1)S_+(u)+R_-(-u-1)R_+(u)}{\psi_-(-u-1)(k+u+1)}\,+\,
\frac{R_+(u)}{k+u+1}\,.
\endmultline
$$
Here the second and the fourth fractions cancel each other, while the
third fraction equals $\psi_+(u)/(k+u+1)$, because of \tht{3.15}.
Consequently, the whole expression is equal to 
$$
-\,\frac{S_-(k)S_+(u)+R_-(k)R_+(u)}{\psi_-(k)(k+u+1)}\,+\,
\frac{\psi_+(u)}{k+u+1}\,.
$$
Substituting this expression instead of the sum in \tht{3.19}, we get
$$
\multline
(WNW)(k,l)=\lim_{u\to l}\left\{
-\,\frac{S_-(k)S_+(u)+R_-(k)R_+(u)}{\psi_-(k)\psi_+(u)(k+u+1)}\,+\,
\frac1{k+u+1}\right\}\\
=-\,\frac{S_-(k)S_+(l)+R_-(k)R_+(l)}{\psi_-(k)\psi_+(l)(k+l+1)}\,+\,
\frac1{k+l+1}\,.
\endmultline
$$
Thus, 
$$
(WNW-W)(k,l)=-\,\frac{S_-(k)S_+(l)+R_-(k)R_+(l)}
{\psi_-(k)\psi_+(l)(k+l+1)}\,
$$
which proves \tht{3.12}. 

This completes the proof of the theorem. \qed
\enddemo

\head \S4 Connection with Meixner polynomials \endhead

In this section we shall show that when one of the parameters $z,\
z'$ becomes an integer, the `++'--block of the hypergeometric kernel
defined in Theorem 3.3 turns into the Christoffel--Darboux kernel for
Meixner orthogonal polynomials.  

The Meixner polynomials form a system $\{\M_n(k;\al+1,\xi)\}$ of
orthogonal polynomials, which corresponds to the following weight
function on $\Z_+$: 
$$
f(k)=f_{\al,\xi}(k)= \frac{(\al+1)_k\xi^k}{k!}
=\frac{\Ga(\al+1+k)\xi^k}{\Ga(\al+1)k!}\,,  \qquad k\in\Z_+\,.
$$
Here $k\in\Z_+$ is the argument and $\al>-1$ and $0<\xi<1$ are
parameters; $\deg \M_n(k;\al+1,\xi)=n$. 
For a detailed information about these polynomials see \cite{NSU}, \cite{KS}.
\footnote{Our normalization of the Meixner polynomials coincides with
that of \cite{KS} and slightly differs from that of \cite{NSU}.}

Meixner polynomials can be expressed through the Gauss hypergeometric
function: 
$$
\gather
\M_n(k;\al+1,\xi)=F(-n,-k;\al+1;\tfrac{\xi-1}{\xi})\\
=\frac{k!\Ga(-\al-n)}{\Ga(1+k-n)\Ga(-\al)}\,
\left( \frac{1-\xi}{\xi}\right)^n\,
F(-n,-\al-n;1+k-n;\w).
\endgather
$$

Basic constants related to these polynomials have the form
$$
\gather
\M_n(k;\al+1,\xi)=a_nk^n\,+\,\{\text{lower degree terms in $k$}\},
\quad a_n=\left( \frac{1-\xi}{\xi}\right)^n\,\frac1{(\al+1)_n}\, \\
h_n=||\M_n(k;\al+1,\xi)||^2=\sum_{k=0}^\infty \M_n^2(k;\al+1,\xi)f(k)
=\frac{n!}{\xi^n(1-\xi)^{\al+1}(\al+1)_n}\,.
\endgather
$$ 

Consider the $N$th Christoffel--Darboux kernel for the Meixner
polynomials. It projects the Hilbert space $\ell^2(\Z_+, f(\cdot))$
on the $N$--dimensional subspace spanned by the polynomials of degree
$\le N-1$. Let us pass from $\ell^2(\Z_+, f(\cdot))$ to the ordinary
$\ell^2$ space on $\Z_+$, which corresponds to the counting measure.
Then the Christoffel--Darboux kernel will be transformed to a certain
kernel, which will be called the {\it Meixner kernel\/} and denoted as
$M_N(k,l)$. We have:
$$
\gather
M_N(k,l)=\sum_{n=0}^{N-1}
\frac{\M_n(k;\al+1,\xi)\,\M_n(l;\al+1,\xi)}{h_n}\,\sqrt{f(k)f(l)} 
=\frac{a_{N-1}}{a_Nh_{N-1}}\, \sqrt{f(k)f(l)}\,\\\times
\frac{\M_N(k;\al+1,\xi)\M_{N-1}(l;\al+1,\xi)
-\M_{N-1}(k;\al+1,\xi)\M_N(l;\al+1,\xi)}{k-l}\,.
\endgather
$$

\proclaim{Proposition 4.1} Let $z=N+\al$, $z'=N$, and let $K_{++}(k,l)$ be
the ``$++$'' block of the corresponding hypergeometric kernel. Then 
$$
K_{++}(k,l)=M_N(k+N,l+N).
$$
\endproclaim

\demo{Proof} The proof is straightforward. \qed  
\enddemo

Consider the $N$--point ``Meixner ensemble'' on $\Z_+$ whose joint
probability distribution has the form 
$$
p(k_1,\dots,k_N)
=const\cdot\prod_{1\le i<j\le N}(k_i-k_j)^2\prod_{i=1}^n f_{\al,\xi}(k_i).
$$ 
The standard argument due to Dyson (see \cite{Dy}, \cite{Me}) shows
that the correlation functions of this ensemble are given by
determinantal formulas with the Meixner kernel: 
$$
\rho_n(x_i,\dots,x_n)=\det[M_N(k_i,k_j)]_{i,j=1}^n.
$$

Then Proposition 4.1 shows that our point process $\Pzxi$ restricted
to the positive copy of $\Z_+$ for $z=N+\al$, $z'=N$ coincides with
the trace of the $N$--point Meixner ensemble on the set
$\{N+1,N+2,\dots\}$. In this subset the number of points of the Meixner
ensemble can vary from 0 to $N$, which agrees with our picture.

\head \S5. Scaling limit of the hypergeometric kernel: the Whittaker
kernel 
\endhead
Recall that the construction of the point processes $\Pzxi$ was
started from certain probability distributions on partitions of an
integer number $n$ denoted as $M_{z,z'}^{(n)}$, see \S1. These
distributions possess an additional important property: they
converge, as $n\to\infty$, to a probability distribution on a certain
limit object $\Omega$ called the {\it Thoma simplex}: 
$$
\Omega=\{\al_1\geq\al_2\geq\ldots\geq 0;\,\be_1\geq\be_2\geq\ldots\geq
0\bigm|\sum_{i=1}^{\infty}(\al_i+\be_i)\leq 1\}\,. \tag5.1
$$    
It is a compact topological space with respect to the
topology of co\-or\-di\-nate--wise convergence.

More precisely, for every $n$ we embed the set $\Y_n$ of partitions
of $n$ into $\Omega$ by making use of the map 
$$
\gathered
\Y_n\ni\la=(p_1,\dots,p_d\,|\,q_1,\dots,q_d)\\
\,\mapsto\,
\left( \frac{p_1+1/2}n\,,
\dots,\frac{ p_d+1/2}n\,,0,0,\dots;\,
\frac{ q_1+1/2}n\,,\dots, 
\frac{q_d+1/2}n\,,0,0,\dots\right)\,.  
\endgathered
\tag5.2
$$
Next, we identify $M_{z,z'}^{(n)}$ with its push--forward under the
map (5.2), so 
that $M_{z,z'}^{(n)}$ turns into a probability measure on $\Omega$
with finite support. 

\proclaim{Theorem 5.1} The measures $M_{z,z'}^{(n)}$ weakly
converge to a probability measure $P_{z,z'}$ on $\Omega$ as $n\to\infty$.
\endproclaim

\demo{Proof} This follows from a general theorem, see
\cite{KOV}. \qed 
\enddemo

Recall now that to construct the process $\Pzxi$  on the
lattice $\Z'$ we have mixed all the distributions $\Mnz$,
$n=0,1,2,\dots$, using the negative binomial distribution with
suitable parameters, see \S1: 
$$
\pi(n)=(1-\xi)^t\,\frac{(t)_n}{n!}\, \xi^n, \qquad \xi\in(0,1).  \tag5.3
$$
Let us  embed $\Z'$ into the punctured line $\R^*=\R\setminus\{0\}$
and then rescale the process $\Pzxi$ by multiplying the coordinates of
its points by $(1-\xi)$. Then the rescaled point configuration in
$\R^*$ that corresponds to $\la\in \Y_n$ 
differs from the image (5.2) of $\la$ in $\Omega$ by the scaling
factor $(1-\xi)n$.  

The discrete distribution on the positive semiaxis with 
$$
\operatorname{Prob}\{(1-\xi)n\}=(1-\xi)^t\,\frac{(t)_n}{n!},\quad
n=0,1,2,\dots 
$$ 
which depends on the parameter $\xi\in(0,1)$, converges, as
$\xi\to 1$, to the gamma--distribution with parameter $t$ 
$$
\gamma(ds)=\frac {s^{t-1}}{\Gamma(t)}\,e^{-s}ds.
\tag 5.4
$$

This brings us to the following construction.
Consider the space $\wt\Omega=\Omega\times \R_+$ with the probability
measure 
$$
\wt P_{z,z'}=P_{z,z'}\otimes \frac {s^{t-1}}{\Gamma(t)}\,e^{-s}ds.
$$
To any point $(\omega=(\alpha|\beta),s)\in\wt\Omega$ we associate a
point configuration in $\R$ as follows 
$$
((\alpha|\beta),s)\mapsto (\al_1 s,\al_2 s,\dots;-\be_1 s,-\be_2 s,\dots).
\tag5.5
$$
Thus, the measure $\wt P_{z,z'}$ defines a point process on $\R^*$
which will be denoted as $\wt{\Cal P}_{z,z'}$ 

Then the considerations above together with Theorem 5.1 suggest the following 

\proclaim{Theorem 5.2}
The point processes $\Pzxi$ scaled by $(1-\xi)$ converge, as $\xi\to
1$, to the point process $\wt{\Cal P}_{z,z'}$. 
\endproclaim

We will give a rigorous formulation of this claim and its
proof in our next paper. 
Meanwhile, we will use this theorem as a prompt. 

The main result of our previous work \cite{P.I} -- \cite{P.V} was an explicit
computation of the correlation functions of $\wt P_{z,z'}$.
\footnote{About correlation functions of point processes living on a
nondiscrete space, see \cite{DVJ}, \cite{L1}, \cite{L2}.}
 To
formulate this result we shall need the classical Whittaker function
$W_{\kappa,\mu}(x)$, $x>0$. 

This function can be characterized as the only solution of the Whittaker
equation
$$
y''-\left(\frac 14 -\frac\kappa x+\frac {\mu^2-\frac 14}{x^2}\right)\, y=0
$$
such that $y\sim x^\kappa e^{-\frac x2}$ as $x\to+\infty$ (see
\cite{E1, Chapter 6}).  Here $\kappa$ and $\mu$ are complex
parameters. Note that 
$$ 
W_{\kappa,\mu}=W_{\kappa,-\mu}.
$$
We shall employ the Whittaker function for real $\kappa$ and real or pure
imaginary $\mu$; then $W_{\kappa,\mu}$ is real.

We introduce the functions
$$
\gathered
\Cal P_\pm(x)=\frac{(zz')^{1/4}}{(\Ga(1\pm z)\Ga(1\pm z')\,x)^{1/2}}\,
W_{\frac{\pm(z+z')+1}2,\frac{z-z'}2}(x), \\
\Cal Q_\pm(x)=\frac{(zz')^{3/4}}{(\Ga(1\pm z)\Ga(1\pm z')\,x)^{1/2}}\,
W_{\frac{\pm(z+z')-1}2,\frac{z-z'}2}(x).
\endgathered   \tag5.6
$$

\proclaim{Theorem 5.3}
The correlation functions of
the process $\wt{\Cal P}_{zz'}$ have the form 
$$ \gathered
\wt{\rho}_n^{(z,z')}(u_1,\ldots,u_n)
=\det\left[ {\Cal K}(u_i,u_j)\right]_{i,j=1}^n,
\\ n=1,2,\ldots;\quad u_1,\ldots,u_n\in \Bbb R^*,
\endgathered
$$
where the kernel ${\Cal K}(u,v)$ is conveniently written in the block form
$$
{\Cal K}(u,v)=\cases
{\Cal K}_{++}(u,v),\quad &u,v>0;\\
{\Cal K}_{+-}(u,-v),\quad &u>0\,,v<0;\\
{\Cal K}_{-+}(-u,v),\quad &u<0\,,v>0;\\
{\Cal K}_{--}(-u,-v),\quad &u,v<0;
\endcases
$$
with
$$
{\Cal K}_{++}(x,y)
=\frac{\Cal P_+(x)\Cal Q_+(y)-\Cal Q_+(x)\Cal P_+(y)}{x-y}\,, $$$$
{\Cal K}_{+-}(x,y)
=\frac{\Cal P_+(x)\Cal P_-(y)+\Cal Q_+(x)\Cal Q_-(y)}{x+y}\,, $$$$
{\Cal K}_{-+}(x,y)
=-\,\frac{\Cal P_-(x)\Cal P_+(y)+\Cal Q_-(x)\Cal Q_+(y)}{x+y}\,,$$$$
{\Cal K}_{--}(x,y)
=\frac{\Cal P_-(x)\Cal Q_-(y)-\Cal Q_-(x)\Cal P_-(y)}{x-y}\,.
$$
\endproclaim
The kernel ${\Cal K}(u,v)$ is called the {\it Whittaker kernel}, see
\cite{P.IV, Th. 2.7}, \cite{BO1, Th. III }. \footnote{In that papers,
the term `Whittaker kernel' concerned the block $\Cal K_{++}$ while the
kernel $\Cal K$ was called the matrix Whittaker kernel. }

Clearly, the hypergeometric kernel (see Theorem 3.3) and the
Whittaker kernel have the same structure. Theorem 5.2 prompts that
the Whittaker kernel is the scaling limit of the hypergeometric one.
In the next theorem we establish this fact by a direct computation. 

\proclaim{Theorem 5.4} For the hypergeometric kernel $K$ given
by Theorem 3.3 and the Whittaker kernel $\Cal K$ given by
Theorem 5.3 the following limit relation holds 
$$
\lim_{\xi\nearrow1}\frac 1{1-\xi}\,
K_{**}\left(\left[\frac u{1-\xi}\right],
\left[\frac v{1-\xi}\right]\right)
=\Cal K_{**}(u,v), \qquad u,v\in \R_+\,,
$$
where the subscript $**$ stands for any of the four symbols $++$,
$+-$, $-+$, $--$. 
\endproclaim 
\demo{Proof}
Take $x,y>0$ and denote 
$$
k= \left[\frac x{1-\xi}\right],\qquad 
l=\left[\frac y{1-\xi}\right].
$$
Then $(1-\xi)k\approx x$, $(1-\xi)l\approx y$.
Since 
$$
\frac1{k-l} \approx \frac{1-\xi}{x-y}\,, \qquad
\frac1{k+l+1} \approx \frac{1-\xi}{x+y}\,,
$$
it is enough to show that
$$
P_\pm(k)\approx \Cal P_\pm(x),\qquad Q_\pm(k)\approx \Cal Q_\pm(x).
\tag 5.7
$$
We shall employ the following asymptotic relation which connects the
hypergeometric function and the Whittaker function: 
$$
\lim_{u\to+\infty}F(a,b;u;1-\tfrac ux)=
x^{\frac{a+b-1}2}e^{\frac x2}W_{\frac{-a-b+1}2,\frac{a-b}2}(x),
\qquad x>0, \tag 5.8
$$
see \cite{E1, 6.8(1)}. Note that $\w=1-\tfrac1{1-\xi}$. Applying
\tht{5.8} we get the following 
limit relations for the hypergeometric functions entering \tht{3.5} and
\tht{3.6}: 
$$
\gather
F(\mp z, \mp z'; k+1; \w)\,\approx\,
x^{\frac{\mp(z+z')-1}2} e^{\frac x2}
W_{\frac{\pm(z+z')+1}2,\frac{z-z'}2}(x), \\
\frac{F(1\mp z, 1\mp z'; k+2; \w)}{(1-\xi)(k+1)}\,\approx\,
x^{\frac{\mp(z+z')-1}2} e^{\frac x2}
W_{\frac{\pm(z+z')-1}2,\frac{z-z'}2}(x).
\endgather
$$
Next, the factor $(\psi_\pm(u))^{1/2}$ entering \tht{3.7} 
behaves as follows ($\psi_\pm(u)$ was defined in (3.1)):
$$
\gathered
(\psi_\pm(k))^{1/2}=\left(t^{1/2}\xi^{k+1/2}\,
\frac{(1\pm z)_k(1\pm z')_k}{k! k!}(1-\xi)^{\pm(z+z')}\right)^{1/2}
\\
\approx\left(\frac{t^{1/2}e^{-x}x^{\pm(z+z')}}
{\Ga(1\pm z)\Ga(1\pm z')}\right)^{1/2}\,.
\endgathered
$$
Finally, from (3.7) we obtain
$$
\align
&P_\pm(k)\approx
\frac{t^{1/4}}{(\Ga(1\pm z)\Ga(1\pm z')x)^{1/2}}\,
W_{\frac{\pm(z+z')+1}2,\frac{z-z'}2}(x)=\Cal P_\pm(x),    \\
&Q_\pm(k)\approx
\frac{t^{3/4}}{(\Ga(1\pm z)\Ga(1\pm z')x)^{1/2}}\,
W_{\frac{\pm(z+z')-1}2,\frac{z-z'}2}(x)=\Cal Q_\pm(x). 
\qed
\endalign
$$
\enddemo
\example{Remark 5.5}
As was demonstrated in \S4, the ``++''--block of the hypergeometric
kernel turns into the Christoffel--Darboux kernel for Meixner
polynomials when $z=N+\alpha$, $z'=N$, $N\in \Z_+$.  
It is well known that in the scaling
limit as $\xi\to 1$, the Meixner polynomials turn into the   
Laguerre polynomials (see \cite{KS}, \cite{NSU}). This agrees
with the fact that for $z=N+\alpha$, $z'=N$, the restriction of the
process $\wt {\Cal P}_{z,z'}$ to the positive semiaxis coincides with the
$N$--point Laguerre ensemble, see \cite{P.III, Remark 2.4}. 

Note that the shift by $N$ which we were doing to match $\Pzxi$ and
the Meixner ensemble disappears after we take the limit. 
\endexample
\example{Remark 5.6} A straightforward check shows that the
scaling limit of the kernel $L(x,y)$ defined by (3.3) is the kernel
$\Cal L(x,y)$ of the operator $\Cal L=\Cal K(1-\Cal K)^{-1}$ where
$\Cal K$ is the integral operator in $L^2(\R^*,dx)$ corresponding to
the Whittaker kernel (the kernel $\Cal L(x,y)$ was explicitly
computed in \cite{P.V, Theorem 2.4}). 
\endexample

\head \S6. Integrable operators
\endhead

In this section we shall show that the operator given by the
Whittaker kernel belongs to the class of integrable operators as
defined by Its, Izergin, Korepin and Slavnov \cite{IIKS}. We shall also
argue that the hypergeometric kernel might be considered as an
example of a {\it discrete} kernel giving an `integrable operator'. 

We shall follow \cite{De} in our description of integrable operators. 

Let $\Sigma$ be an oriented contour in $\C$. We call an operator $V$
acting in $L^2(\Sigma,|d\zeta|)$ {\it integrable} if its kernel has
the form 
$$
V(\zeta,\zeta')
=\frac {\sum_{j=1}^N f_j(\zeta)g_j(\zeta')}{\zeta-\zeta'},
\quad \zeta,\ \zeta'\in\Sigma, 
$$
for some functions $f_j,\ g_j$, $j=1,\dots,N$. We shall always assume
that $$\sum_{j=1}^N f_j(\zeta)g_j(\zeta)=0,\quad \zeta\in\Sigma,$$ so
that the kernel $V(\zeta,\zeta')$ is nonsingular (this assumption is
not necessary for the general theory).

The notion of an integrable operator was first introduced in \cite{IIKS}. 

It turns out that for an integrable operator $V$ the operator
$R=V(1+V)^{-1}$ is also integrable.
 
\proclaim{Proposition 6.1 \cite{IIKS}} Let $V$ be an integrable operator
as described above and $R=1-(1+V)^{-1}=V(1+V)^{-1}$. Then the kernel
$R(\zeta,\zeta')$ has the form 
$$
R(\zeta,\zeta')
=\frac {\sum_{j=1}^N F_j(\zeta)G_j(\zeta')}{\zeta-\zeta'},
\quad \zeta,\zeta'\in\Sigma, 
$$
where
$$
F_j=(1+V)^{-1}f_j,\qquad G_j=(1+V^t)^{-1}g_j,\quad j=1,\dots,N.
$$
If $\sum_{j=1}^N f_j(\zeta)g_j(\zeta)=0$ on $\Sigma$, then 
$\sum_{j=1}^N F_j(\zeta)G_j(\zeta)=0$ on $\Sigma$ as well.
\endproclaim

\demo{Proof} See \cite{KBI, ch. XIV}, \cite{De}. \qed 
\enddemo

It is not difficult to show that for integrable operators $V_1$,
$V_2$, the product $V_1V_2$ is also integrable. This fact and
Proposition 6.1 imply that operators of the form $I+V$ where $V$ is
integrable form a group.  

A remarkable fact is that the function $F_j$, $G_j$ can be expressed
via a suitable Riemann--Hilbert problem, see \cite{IIKS}, \cite{De}
for details.  

Now we pass to a much more special situation. Let $\Sigma=\R^*$.
According to the splitting $\R^*=\R_+\sqcup \R_-$ and further
identification of $\R_-$ with a second copy of $\R_+$, we shall
sometimes write the kernels of operators in $L^2(\R^*)$ in block
form. 

Consider an integral operator $V$ on $\R^*$ whose kernel $V(x,y)$ has
the following block form
$$
[V](x,y)=\left[\matrix 0&\frac {h_{+}(x)h_-(y)}{x+y}\\
-\frac {h_{+}(y)h_-(x)}{x+y}&0\endmatrix\right]\,, 
\qquad  x>0, y>0,   \tag 6.1
$$
for some functions $h_+(x)$ and $h_-(x)$ defined on the positive semiaxis.
Then the operator $V$ is integrable. Indeed,
$$
V(x,y)=\frac{f_1(x)g_1(y)+f_2(x)g_2(y)}{x-y},\qquad x,y\in \R^*
$$
where
$$
f_1(x)=\cases 0,&x>0\\h_-(-x),&x<0\endcases,\quad
f_2(x)=\cases h_+(x),&x>0\\0,&x<0\endcases
$$$$
g_1(x)=\cases h_+(x),&x>0\\0,&x<0\endcases,\quad
g_2(x)=\cases 0,&x>0\\h_-(-x),&x<0\endcases
$$

Assume that there exist four functions $A_\pm(x)$, $B_\pm(x)$ defined
on the positive semiaxis such that 
$$
\widehat A_\mp=\frac{B_\pm}{h_\pm^2},\quad 
\widehat B_\mp=1-\frac{A_\pm}{h_\pm^2}
\tag 6.3
$$
where
$$
\widehat\varphi(x)=\int_{y>0}\frac {\varphi(y)dy}{x+y}
\tag 6.4
$$
is the Stieltjes transform.

Then Proposition 6.1 implies that the kernel of the operator
$R=V(1+V)^{-1}$ has the form 
$$
R(x,y)=\frac{F_1(x)G_1(y)+F_2(x)G_2(y)}{x-y},\qquad x,y\in \R^*
$$
with
$$
F_1(x)=\cases-\frac {B_+(x)}{h_+(x)} ,&x>0\\  
\frac {A_-(-x)}{h_-(-x)},&x<0\endcases;\qquad
F_2(x)=\cases\frac {A_+(x)}{h_+(x)} ,&x>0\\ 
\frac{B_-(-x)}{h_-(-x)} ,&x<0\endcases;
$$
$$
G_1(x)=\cases \frac {A_+(x)}{h_+(x)},&x>0\\  
-\frac{B_-(-x)}{h_-(-x)},&x<0\endcases;\qquad
G_2(x)=\cases \frac {B_+(x)}{h_+(x)} ,&x>0\\  
\frac{A_-(-x)}{h_-(-x)},&x<0\endcases.
$$

In block form the kernel $R(x,y)$ can be written as follows:
$$
[R](x,y)=\left[\matrix
\frac 1{h_+(x)h_+(y)}\frac{A_+(x)B_+(y)-B_+(x)A_+(y)}{x-y}&
\frac 1{h_+(x)h_-(y)}\frac{A_+(x)A_-(y)+B_+(x)B_-(y)}{x+y}\\
\frac 1{h_-(x)h_+(y)}\frac{-A_-(x)A_+(y)-B_-(x)B_+(y)}{x+y}&
\frac 1{h_-(x)h_-(y)}\frac{A_-(x)B_-(y)-B_-(x)A_-(y)}{x-y}
\endmatrix\right] 
$$

All these formulas work perfectly well for the Whittaker kernel. If
we set, using the notation of \S5, 
$$
h_\pm(x)=\frac {\sqrt{\sin\pi z\,\sin\pi z'}}
{\pi}\,x^{\mp \frac {z+z'}2}e^{-x/2};
$$
$$
A_\pm(x)={h_{\pm}(x)}{\Cal P_\pm(x)},\quad   
B_\pm(x)={h_{\pm}(x)}{\Cal Q_\pm(x)},
$$
then the kernel $R(x,y)$ coincides with the Whittaker kernel $\Cal
K(x,y)$. The form (6.1) of the kernel of $V=\Cal K(1-\Cal K)^{-1}$
was obtained in \cite{P.V}. 
The formulas (6.3) in this case can be derived from the known
formulas for the Stieltjes transform of the (suitably normalized)
Whittaker function \cite{E2, 14.3(53)}. 

It is a remarkable fact that all the formulas above also work for the
hypergeometric kernel. This kernel lives on the lattice
$\Z'\times\Z'$, so one can call the operator corresponding to the
hypergeometric kernel a {\it discrete integrable operator}.  

Indeed, exact expressions for $h_+(k)$ and $h_-(k)$ can be easily extracted
from (3.3) if we take $V=L$: 
$$
h_\pm(k)=(\psi_\pm(k))^{1/2}. 
$$
Then, as before, we have, cf. (3.5)--(3.7),
$$
 A_\pm(k)={h_{\pm}(k)}{P_\pm(k)}=R_\pm(k),\quad  
B_\pm(k)={h_{\pm}(k)}{Q_\pm(k)}= S_\pm(k),
$$
and the kernel $R(x,y)$ coincides with the hypergeometric kernel $K(x,y)$.

Relations (6.3) are exactly the relations of Lemma 3.4, see (3.14).

If we consider the continuous case, then from the general theory of
Riemann--Hilbert problems one can extract the following identity for 
the analytic continuations of $B_\pm,\ A_\pm$:
$$
A_+(\zeta)A_-(-\zeta)+B_+(\zeta)B_-(-\zeta)
=h^2_+(\zeta)h_-^2(-\zeta),\quad \zeta\in \C\setminus\R.
\tag 6.5
$$

It means that the determinant of the solution of the corresponding
Riemann--Hilbert problem is identically equal to 1 (this follows from
the fact that the determinant of the corresponding jump matrix is
identically equal to 1). 

Though we do not have an analog of the Riemann--Hilbert problem in
the discrete case, a discrete analog of (6.5) still holds, see Lemma
3.5. 

\example{Remark 6.2} Both the Whittaker and the hypergeometric kernel
possess the symmetry property 
$$
R(x,y)=\operatorname{sgn}(x)\operatorname{sgn}(y)R(y,x)
\tag 6.6
$$
which, perhaps, emerged for the first time: most (integrable) kernels
arising in Random Matrix Theory and mathematical physics are simply
symmetric. 
The formula (6.6) means that the corresponding integral operator in
$L^2(\R,dx)$ or in $\ell^2(\Z')$ is symmetric with respect to an
indefinite inner product.  
\endexample

\head \S7. Appendix: some relations for the Gauss hypergeometric function
\endhead
In this section we prove Lemmas 3.4 and 3.5. First, we shall
reformulate them using (3.5) and (3.6).
\proclaim{Lemma 3.4'} For $\xi\in(0,1)$ the following
decompositions hold: 
$$
\gather
\frac{F(a,b;u+1;\w)}u=
\sum_{k=0}^\infty 
\frac{(a)_k (b)_k\xi^k(1-\xi)^{a+b-1}}{k!\,k!\,(u+k)}\,
F(1-a,1-b;k+1;\w),\\
1-F(a,b;u;\w)=
\sum_{k=0}^\infty 
\frac{(a)_{k+1}(b)_{k+1}\xi^{k+1}(1-\xi)^{a+b-1}}
{k!\,k!\,(u+k)}\,
\frac{F(1-a,1-b;k+2;\w)}{k+1}. 
\endgather
$$
Specifically, the series in the RHS absolutely converges for
$u\ne0,-1,-2,\dots$ and represents a meromorphic function; the
both formulas are viewed as equalities of meromorphic functions in
$u$. 
\endproclaim
\proclaim{Lemma 3.5'}
$$
\gathered
F(-z,-z';u+1;w)F(z,z';-u;w)\\
+zz'w(1-w)\,\frac{F(-z+1,-z'+1;u+2;w)}{u+1}\,
\frac{F(z+1,z'+1;-u+1;w)}{u}\,=\,1. 
\endgathered 
\tag 7.1
$$
\endproclaim
\demo{Proof of Lemma 3.4'} Let us check the first relation. The RHS
has the form 
$$
\sum_{k=0}^\infty \frac{A_k}{u+k},
$$
where the coefficients $A_k$ rapidly decrease as $k\to\infty$,
because of the factor $\xi^k$ (the factor $F(1-a,1-b;k+2;\w)$ remains
bounded as $k\to+\infty$, and the remaining expression has at most
polynomial growth in $k$). Consequently, the RHS is indeed a 
converging series representing a meromorphic function in $u$. This function
has simple poles at $u=0,-1,-2, \dots$. Using the formula
$$
\gathered
\Res_{c=-k} F(a,b;c;w)=
\frac{(a)_{k+1}(b)_{k+1}(-1)^kw^{k+1}}{(k+1)!k!}
F(a+k+1,b+k+1;k+2;w)\\
=\frac{(a)_{k+1}(b)_{k+1}(-1)^kw^{k+1}(1-w)^{-a-b-k}}{(k+1)!k!}
F(1-a,1-b;k+2;w) 
\endgathered
\tag 7.2
$$
for the residues of the hypergeometric function,
 one verifies that
the residues of the RHS at $u=-1,-2,\dots$ are the same as for the LHS, and it
is directly seen that the residues at $u=0$ coincide, too. Moreover,
the same claim holds not only for $\xi\in(0,1)$ but for any complex
$\xi$ ranging over the unit disc $|\xi|<1$, and the both sides are
holomorphic in $\xi$. 

Let us expand both the LHS and the RHS into Taylor series
in $\xi$ and compare the respective Taylor coefficients.
Each Taylor coefficient (on the left and on the right), viewed as a
function in $a,b,u$, is a rational expression which is polynomial in
$a,b$. This implies that it suffices to prove our relation,
say, for $a=0,-1,-2,\dots$. 

Thus, we may assume that $a=0,-1,-2,\dots$. For these values of $a$,
the LHS becomes a rational function in $u$, and the series in the RHS
terminates and, consequently, is a rational function in $u$, too.
Next, we know that the both sides have the same singularities.
Finally, they both behave as $O(1/u)$ as $|u|\to\infty$ 
(indeed, as was mentioned above, the hypergeometric function in the
numerator of the LHS is $1+O(1/u)$, so that the whole expression is
$O(1/u)$, and for the RHS the same holds, because the series
terminates).  Consequently, the both sides are identical. This
concludes the proof of our first relation.  

The second relation is verified similarly. \qed
\enddemo
\demo{Proof of Lemma 3.5'}
We use the same argument as in the previous proof. The desired
relation can be viewed as an equality of 
power series in the variable $w$. The coefficients of the series are
polynomials in $z,z'$, so that we may assume, without loss of
generality, that one of the parameters $z,z'$ takes positive integer
values while another takes negative integer values. In this case all
the four hypergeometric series entering our relation terminate and,
so, are rational functions in $u$.  

Next, we examine the possible singularities of the LHS of \tht{7.1}.
Here only simple poles at the points $u\in\Z$ may occur, but it turns
out that the residue at any $u\in\Z$ vanishes. That is, the
contributions of the two products cancel each other. To see this, one
can use, for example, \tht{7.2}. 

Finally, we remark that, under our specialization of $z,z'$, the LHS of
\tht{7.1} is $1+O(1/u)$ as $|u|\to\infty$ . 
This concludes the proof.  \qed
\enddemo

\bigskip
\bigskip

{\smc A.~Borodin}: Department of Mathematics, The University of
Pennsylvania, Philadelphia, PA 19104-6395, U.S.A.  

E-mail address:
{\tt borodine\@math.upenn.edu}

{\smc G.~Olshanski}: Dobrushin Mathematics Laboratory, Institute for
Problems of Information Transmission, Bolshoy Karetny 19, 101447
Moscow GSP-4, RUSSIA.  

E-mail address: {\tt olsh\@iitp.ru,
olsh\@glasnet.ru}


\Refs 
\widestnumber\key{BDJ1}

\ref\key BDJ1
\by J.~Baik, P.~Deift, K.~Johansson
\paper On the distribution of the length of the longest increasing
subsequence of random permutations
\paperinfo Preprint, 1998, available via 
{\tt http://xxx.lanl.gov/ abs/math/9810105}
\endref

\ref\key BDJ2
\bysame
\paper On the distribution of the length of the second row of a Young
diagram under Plancherel measure
\paperinfo Preprint, 1999, available via 
{\tt http://xxx.lanl.gov/abs/math/ 9901118}
\endref

\ref\key BOO
\by A.~Borodin, A.~Okounkov and G.~Olshanski
\paper Paper in preparation
\endref

\ref\key BO1
\by A.~Borodin and G.~Olshanski
\paper Point processes and the infinite symmetric group 
\jour Math. Research Lett.
\vol 5
\yr 1998
\pages 799--816 (preprint version available via 
{\tt http://xxx.lanl.gov/abs/ math/9810015}) 
\endref

\ref\key BO2
\by A.~Borodin and G.~Olshanski
\paper Paper in preparation
\endref

\ref\key DVJ
\by D.~J.~Daley, D.~Vere--Jones 
\book An introduction to the theory of point processes 
\bookinfo Springer series in statistics 
\publ Springer 
\yr 1988 
\endref 

\ref\key De
\by P.~Deift
\paper Integrable operators
\paperinfo Preprint, February 23, 1998
\endref

\ref\key Dy
\by F.~J.~Dyson
\paper  Statistical theory of the energy levels of complex systems I, II, III
\jour J. Math. Phys. 
\vol 3
\yr 1962
\pages 140-156, 157-165, 166-175
\endref

\ref\key E1
\by A.~Erdelyi (ed.) 
\book Higher transcendental functions, {\rm Vol. 1}
\publ Mc Graw--Hill
\yr 1953
\endref

\ref\key E2
\by A.~Erdelyi (ed.) 
\book Tables of integral transforms, {\rm Vol. 1}
\publ McGraw--Hill
\yr 1954
\endref

\ref\key IIKS
\by A.~R.~Its, A.~G.~Izergin, V.~E.~Korepin, N.~A.~Slavnov
\paper Differential equations for quantum correlation functions
\jour Intern. J. Mod. Phys. 
\vol B4
\yr 1990
\pages 1003--1037
\endref

\ref\key J
\by K.~Johansson
\paper Shape fluctuations and random matrices
\paperinfo Preprint, 1999, available via 
{\tt http://xxx.lanl.gov/abs/math/9903134}
\endref

\ref\key K
\by S.~V.~Kerov
\paper Anisotropic Young diagrams and Jack symmetric functions
\jour Funct. Anal. Appl.
\pages to appear (preprint version available via 
{\tt http:/xxx.lanl.gov/abs/math/9712267})
\endref

\ref\key KOO
\by S.~Kerov, A.~Okounkov, G.~Olshanski
\paper The boundary of Young graph with Jack edge multiplicities
\jour Intern. Math. Res. Notices 
\yr 1998
\issue 4
\pages 173--199 
\endref

\ref \key KOV 
\by S.~Kerov, G.~Olshanski, A.~Vershik 
\paper Harmonic analysis on the infinite symmetric group. A deformation 
of the regular representation 
\jour Comptes Rend. Acad. Sci. Paris, S\'er. I 
\vol 316 
\yr 1993 
\pages 773-778 
\endref 

\ref\key KV 
\by S.~V.~Kerov and A.~M.~Vershik 
\paper The Grothendieck group of the infinite symmetric group and 
symmetric functions with the elements of the $K_0$-functor theory 
of AF-algebras 
\inbook Representation of Lie groups and related topics 
\bookinfo Adv. Stud. Contemp. Math. {\bf 7} 
\eds A.~M.~Vershik and D.~P.~Zhelobenko 
\publ Gordon and Breach 
\yr 1990 
\pages 36--114
\endref

\ref\key KS
\by R.~Koekoek and R.~F.~Swarttouw
\paper The Askey--scheme of hypergeometric orthogonal polynomials
and its $q$-analogue
\paperinfo available via  
{\tt ftp://ftp.twi.tudelft.nl/TWI/publications/ tech-reports/1998/
DUT-TWI-98-17.ps.gz}
\endref

\ref \key KBI 
\by  V.~E.~Korepin, N.~M.~Bogoliubov, A.~G.~Izergin
\book Quantum inverse scattering method and correlation functions
\publ Cambridge University Press
\yr 1993
\endref

\ref\key L1
\by A.~Lenard
\paper Correlation functions and the uniqueness of the state in classical 
statistical mechanics
\jour Comm. Math.Phys
\vol 30 
\yr 1973
\pages 35--44
\endref

\ref\key L2
\by A.~Lenard
\paper States of classical statistical mechanical systems of
infinitely many paticles. I, II
\jour Archive for Rational Mech. Anal. 
\vol 59
\yr 1975
\pages 219--239, 241--256
\endref 
 

\ref\key Ma 
\by I.~G.~Macdonald 
\book Symmetric functions and Hall polynomials 
\bookinfo 2nd edition 
\publ Oxford University Press 
\yr 1995 
\endref 

\ref \key Me 
\by M.~L.~Mehta
\book Random matrices
\publ 2nd edition, Academic Press, New York
\yr 1991
\endref

\ref\key NSU
\by A.~F.~Nikiforov, S.~K.~Suslov and V.~B.~Uvarov
\book Classical orthogonal polynomials of a discrete variable
\bookinfo Springer Series in Computational Physics
\publ Springer
\publaddr New York
\yr1991
\endref

\ref\key O
\by A.~Okounkov
\paper Random matrices and random permutations
\paperinfo Preprint, 1999, available via 
{\tt http://xxx.lanl.gov/abs/math/9903176}
\endref

\ref\key  P.I
\by G.~Olshanski
\paper Point processes and the infinite symmetric group. Part I: The
general formalism and the density function
\paperinfo Preprint, 1998, available via 
{\tt http://xxx.lanl.gov/abs/ math/9804086}
\endref

\ref\key P.II
\by A.~Borodin
\paper Point processes and the infinite symmetric group. Part II:
Higher correlation functions
\paperinfo Preprint,
 1998, available via
 {\tt http://xxx.lanl.gov/abs/math/9804087}
\endref

\ref\key P.III
\by A.~Borodin and G.~Olshanski
\paper Point processes and the infinite symmetric group. Part III:
Fermion point processes
\paperinfo Preprint, 1998, available via 
{\tt http://xxx.lanl.gov/abs/math/ 9804088}
\endref

\ref\key P.IV
\by A.~Borodin
\paper Point processes and the infinite symmetric group. Part IV:
Matrix Whittaker kernel
\paperinfo Preprint,
 1998, available via  
{\tt http://xxx.lanl.gov/abs/math/9810013}
\endref

\ref\key P.V
\by G.~Olshanski
\paper Point processes and the infinite symmetric group. Part V:
Analysis of the matrix Whittaker kernel 
\paperinfo Preprint, 1998, available via 
{\tt http://xxx.lanl.gov/abs/math/ 9810014}
\endref


\ref\key T
\by E.~Thoma
\paper Die unzerlegbaren, positive--definiten Klassenfunktionen
der abz\"ahlbar unendlichen, symmetrischen Gruppe
\jour Math.~Zeitschr.
\vol 85
\yr 1964
\pages 40-61
\endref

\ref\key TW
\by C.~A.~Tracy and H.~Widom
\paper Level spacing distributions and the Airy kernel
\jour Comm. Math. Phys.
\vol 159
\yr 1994
\pages 151--174
\endref

\ref\key V
\by A.~M.~Vershik
\paper Statistical mechanics of combinatorial partitions, and their
limit shapes
\jour Funct. Anal. Appl. 
\vol 30
\yr 1996
\pages 90--105
\endref

\ref\key VK
\by A.~M.~Vershik, S.~V.~Kerov
\paper Asymptotic theory of characters of the symmetric group
\jour Funct. Anal. Appl. 
\vol 15
\yr 1981
\pages 246--255
\endref

\endRefs
\enddocument
\end